\documentclass[12pt,reqno]{amsart}
\usepackage{xy, amssymb, amsfonts, amsbsy, amsthm, amsmath, amscd, latexsym, stmaryrd, epic, eepic,eucal}

\DeclareMathAlphabet{\mathpzc}{OT1}{pzc}{m}{it}

\newenvironment{dem}{\begin{proof}[\bf Proof]}{\end{proof}}

\newtheorem{theorem}{\bf Theorem}[section]
\newtheorem{lemma}[theorem]{\bf Lemma}
\newtheorem{propos}[theorem]{\bf Proposition}

\theoremstyle{definition}
\newtheorem{defi}[theorem]{\bf Definition}
\newtheorem{oss}[theorem]{\bf Remark}

\newtheorem{exm}[theorem]{\bf Example}

\newcommand{\B}{\text{B}}

\newcommand{\C}{\mathbb C}
\newcommand{\Cg}{\mathfrak C}

\newcommand{\Gm}{\mathbb G_{\textbf{m}}}
\newcommand{\Ga}{\mathbb G_{\textbf{a}}}

\newcommand{\M}{\mathfrak M}

\newcommand{\N}{\mathbb N}
\newcommand{\Nf}{\mathcal N}
\newcommand{\Of}{\mathcal O}
\newcommand{\Pro}{\mathbb P}
\newcommand{\Q}{\mathbb Q}

\newcommand{\St}{\mathcal S}

\newcommand{\aut}{\text{Aut}}

\newcommand{\cart}{\ar @{} [dr] |{\Box}}
\newcommand{\cat}{\mathpzc{Sch}_{\mathbb C}}

\newcommand{\cu}{\mathpzc C}

\newcommand{\id}{\text{id}}

\newcommand{\km}{\mathpzc k}

\newcommand{\spc}{\text{Spec} \mathbb C}
\newcommand{\spe}{\text{Spec}}
\newcommand{\spk}{\text{Spec}\Omega}

\newcommand{\unoa}{\xymatrix { *=0{\bullet} \ar@{-}[r] & *=0{\bullet} }}
\newcommand{\duea}{\xymatrix { *=0{\bullet} \ar@{-}[r] & *=0{\bullet} \ar@{-}[r] &*=0{\bullet} }}

\newcommand{\trea}{\xymatrix { *=0{\bullet} \ar@{-}[r] & *=0{\bullet} \ar@{-}[r] &*=0{\bullet} \ar@{-}[r] &*=0{\bullet} }}

\newcommand{\treb} {\xymatrix @R=8pt { & *=0{\bullet} \ar@{-}[dd] & \\ & &  \\ & *=0{\bullet} \ar@{-}[dl] \ar@{-}[dr] & \\ *=0{\bullet} & & *=0{\bullet}}}

\newcommand{\trebg}{\xymatrix @R=8pt { & *=0{} \ar@{-}[4,0] & \\ *=0{} \ar@{-}[0,2] & & *=0{} \\ *=0{} \ar@{-}[0,2] & & *=0{}\\ *=0{} \ar@{-}[0,2] & & *=0{}\\ & *=0{} &}}

\newcommand{\quattroa}{\xymatrix { *=0{\bullet} \ar@{-}[r] & *=0{\bullet} \ar@{-}[r] &*=0{\bullet} \ar@{-}[r] &*=0{\bullet} \ar@{-}[r] &*=0{\bullet} }}

\newcommand{\quattrob}{\xymatrix {& *=0{\bullet}  \ar@{-}[d]& \\ *=0{\bullet} \ar@{-}[r] & *=0{\bullet} \ar@{-}[r] \ar@{-}[d] & *=0{\bullet}\\ & *=0{\bullet} &}}

\newcommand{\quattroc}{\xymatrix @R=4pt {*=0{\bullet} \ar@{-}[dr] & & &\\ & *=0{\bullet} \ar@{-}[r] &*=0{\bullet} \ar@{-}[r] &*=0{\bullet}\\ *=0{\bullet} \ar@{-}[ur] & & &}}

\newcommand{\quattrocg}{\xymatrix @R=8pt { & *=0{} \ar@{-}[4,0] & & &\\ *=0{} \ar@{-}[0,2] & & *=0{} &*=0{} \ar@{-}[2,0] & \\ *=0{} \ar@{-}[0,4] & & & &*=0{}\\ *=0{} \ar@{-}[0,2] & & *=0{} &*=0{} &\\ & *=0{} & & &}}

\newcommand{\cinquenz}{\xymatrix { *=0{\bullet} \ar@{-}[d] & *=0{\bullet} \ar@{-}[d] \\ *=0{\bullet} \ar@{-}[r] & *=0{\bullet} \\
*=0{\bullet} \ar@{-}[u] & *=0{\bullet} \ar@{-}[u]}}

\input xy
\xyoption{all}

\begin{document}

\title[The Stack of Rational Curves]{THE STACK OF RATIONAL NODAL CURVES}
\author{Damiano Fulghesu}
\address{Department of Mathematics,
University of Missouri,
Columbia, MO 65211}
\email{damiano@math.missouri.edu} 

\begin{abstract}
In this series of three papers we start to investigate the rational Chow ring of the stack $\M_{0}$ consisting of nodal curves of genus $0$, in particular we determine completely the rational Chow ring of 
the substack $\M_{0}^{\leq 3}$ consisting of curves with at most $3$ nodes. In the first paper we construct the stack $\M_0$ and its stratification by nodes and show that the map from the universal curve to the stack is not representable in the category of schemes. 
\end{abstract}

\maketitle

\section{Introduction}

\medskip

Foundations of intersection theory on moduli spaces of curves were laid in the famous Mumford's paper \cite{mumenum} where he posed the problem for stable curves and carried out computations for genus 2 curves.

Afterward, in a series of papers,  Faber \cite{Fab} computed $A^*_{\Q} \overline{\mathcal M}_3$ and gave some partial results on $A^*_{\Q} \overline{\mathcal M}_4$. Subsequently Izadi \cite{Iza} proved that $A^*\overline{\mathcal M}_5$  is generated by tautological classes.

The theory is more subtle for spaces of very low genus such as 0 and 1 for which the automorphism group of smooth curves is not finite. One possible approach is to consider the space of stable $n$-pointed curves. Within this context Keel \cite{Keel} exhibits a complete description of the integral Chow ring of $\overline{\mathcal M}_{0,n}$.

In a series of papers we inquire directly the moduli stack of (not pointed) rational curves $\M_0$ using equivariant intersection theory for quotient stacks as developed by Edidin and Graham \cite{EGRR} and Kresch's work \cite{kre} on general Artin stacks.

We make use of the stack structure since the associated moduli space of an algebraic stack loses much information. For example the stack $\M^0_0$ of rational smooth curves has a moduli space which consists of a point (all smooth rational curves are isomorphic) but as a stack it is the classifying space of the automorphism group of $\Pro^1$ (i.e. $\B \Pro \text{GL}_2$) which has a much deeper  structure.

The purpose of this paper is to prove that $\M_0$ can be stratified into smooth quotient substacks $\M_0^{\Gamma}$ defining families of curves which are topologically equivalent. This allows us to apply Kresch's theory and define a Chow ring for $\M_0$.

Unfortunately even defining tautological classes is extremely difficult since the morphism from the universal curve $\mathcal C \to \M_0$ is not projective (not even representable in the category of schemes as shown in Example \ref{HV}) and we lack a sufficiently general definition for pushforward of Chern classes. Moreover it is still an open question whether, for $n \geq 2$, the stack $\M^{\leq n}_0$ of rational curves with at most $n$ nodes is a quotient stack. In the second paper \cite{fulg2}, we explain how to use Grothendieck-Riemann-Roch to define tautological Mumford classes $\kappa_i$  on $\M_0^{\leq 3}$.

In the third paper \cite{fulg3} we give a presentation for $A^*(\M_0^{\leq 3})\otimes \Q$. There are 10 generators: the classes corresponding to the five topological types of curves with at most 3 nodes,  the Mumford class $\km_{2}$ and $4$ other generators which are somewhat new with respect to the tautological classes introduced for stable curves. 

For completeness we mention that the problem of computing integral Chow rings of rational curves with at most $n$ nodes is solved for the lowest cases $n=0$ \cite{Pa} and $n=1$ \cite{EF}. It is still unknown for larger $n$.

{\bf Acknowledgments.} I am grateful to my thesis advisor Angelo Vistoli for his patient and constant guidance and Rahul Pandharipande who suggested the problem to him.

I also wish to thank Dan Edidin for very helpful remarks.

\section{Description of the stack $\M_0$}

Let $\cat$ be the site of schemes over the complex point $\spc$ equipped with the \'etale topology.

\begin{defi}
We define the category of rational nodal curves $\M_0$ as the category over $\cat$ whose objects are flat and proper morphisms of finite presentation $C \xrightarrow{\pi} T$ (where $C$ is an algebraic space over $\C$ and $T$ is an object in $\cat$) such that for every geometric point $\spk \xrightarrow{t} T$ (where $\Omega$ is an algebraically closed field) the fiber $C_{t}$
\begin{equation*}
\xymatrix{
C_{t} \cart \ar[r]^{t'} \ar[d]_{\pi'} &C \ar[d]^{\pi}\\
\spk \ar[r]^t &T
}
\end{equation*}
is a projective reduced nodal curve such that $h^0(C_t, \Of_{C_t}) = 1$ and $h^1(C_t, \Of_{C_t}) = 0$. Morphisms in $\M_0$ are cartesian diagrams
\begin{equation*}
\xymatrix{
C' \cart \ar[r] \ar[d]_{\pi'} &C \ar[d]^{\pi}\\
T' \ar[r] &T
}
\end{equation*}
The projection $pr : \M_0 \xrightarrow{} \cat$ is the forgetful functor
\begin{eqnarray*}
pr:
&
 \left \{
\vcenter {\xymatrix{
C' \cart  \ar[r] \ar[d]_{\pi'} & C \ar[d]^{\pi}\\
T' \ar[r] & T
}}
\right \}
&
\mapsto \{ T' \xrightarrow{} T \}
\end{eqnarray*}
\end{defi}

We refer to \cite[Proposition 1.10]{fulg1} for the proof that $\M_0$ is an Artin stack  in the sense of
\cite{Artd} (the proof uses standard arguments).

Here we describe more precisely the geometric points of $\M_0$. Given an algebraically closed field $\Omega$, let us define a {\it rational tree (on $\Omega$)} to be a connected nodal curve with finite components each of them is isomorphic to $\Pro^1_\Omega$ and which has no closed chains.
Classical cohomology arguments allow us to state that geometric points of $\M_0$ are rational trees.

\begin{oss}It is important to notice that by definition geometric fibers of a family of curves $C \xrightarrow{\pi}T$ in $\M_0$ are projective and therefore are schemes. Moreover it is further known (see \cite {Kn} V Theorem 4.9) that  a curve (as an algebraic space) over an algebraically closed field $\Omega$ is a scheme.

Notwithstanding the fact that we consider families of curves whose fibers and bases are schemes we still need algebraic spaces because, as we show below, there exist families of rational nodal curves $C \xrightarrow{\pi} T$ where $C$ is not a scheme.
\end{oss}

\begin{exm} \label{HV}
Let us consider the following example which is based on Hironaka's example \cite{Hir} of an analytic threefold which is not a scheme (this example can be found also in \cite{Kn}, \cite{hrt} and \cite{Shf}).

Let $S$ be a projective surface over $\C$ and $i$ an involution without fixed points. Suppose that there is a smooth curve $C$ such that $C$ meets transversally the curve $i(C)$ in two points $P$ and $Q$. Let $M$ be the product $S \times \Pro ^1_{\C}$ and let $e: S \to M$ be the embedding $S \times 0$. On $M - e(Q)$, first blow up the curve $e(C - Q)$ and then blow up the strict transform of $e(i(C) - Q)$. On $M - e(P)$, first blow up the curve $e(i(C) - P)$ and then blow up the strict transform of $e(C - P)$. We can glue these two blown-up varieties along the inverse images of $M - e(P) - e(Q)$. The result is a nonsingular complete scheme $\widetilde M$. 

We have  an action of $\Cg_2$ (the order-2 group generated by the involution $i$) on $M$ induced by the action of $\Cg_2$ on $S$ and the trivial action on $\Pro^1_{\C}$. Furthermore we can lift the action of $\Cg_2$ to $\widetilde M$ such that we obtain an equivariant map
$$
\widetilde M \xrightarrow{\pi} S;
$$
we still call $i$ the induced involution on $\widetilde M$. We have an action of $\Cg$ on $\widetilde M$ which is faithful and we obtain geometric quotients $f$ and $g$
\begin{equation*}
\xymatrix{
\widetilde M \ar[d]_{\pi} \ar[r]^{f} & \widetilde M /\Cg_2  \ar[d]^{\widetilde{\pi}} \\
S \ar[r]^{g} & S/\Cg_2
}
\end{equation*}

Now  we notice that $\pi$ is a family of rational nodal curves. In particular the generic geometric fiber is isomorphic to $\Pro^1_{\C}$. Geometric fibers on the two curves $C$ and $i(C)$ have one node except those on $P$ and $Q$ which have two nodes. Similarly we have that also $\widetilde{\pi}$ is a family of rational nodal curves whose fibers with a node are on $g(C)=g(i(C))$ and there is only one fiber on $g(P)=g(Q)$.

Since $\Cg_2$ acts faithfully on $S$ and $S$ is projective we have that $S/\Cg_2$ is a scheme. Furthermore  we mention (see \cite{Kn} Chapter 4) that the category of separated algebraic spaces  is stable under finite group actions. Consequently the family $\widetilde{\pi}$ belongs to $\M_0(S/ \Cg_2)$. But $\widetilde M / \Cg_2$ is not a scheme (to prove this we can follow the same argument of \cite{hrt} Appendix B Example 3.4.2).

In order to complete our example, we have to exhibit a surface which satisfies the required properties. Take the Jacobian $J_2=\text{Div}^0(C)$ of a genus 2 curve $C$. 
Let us choose on $C$ two different points $p_0$ and $q_0$ such that $2(p_0 - q_0) \sim 0$. Let us consider the embedding
\begin{eqnarray*}
C &\to& J_2\\
p \in C &\mapsto& |p-q_0|.
\end{eqnarray*}
We still call $C$ the image of the embedding. Let $i$ be the translation on $J_2$ of $|p_0 - q_0|$. By definition $i$ is an involution such that $C$ and $i(C)$ meets transversally in two points: $0$ and $|p_0 - q_0|$.
\end{exm}

\section{Stratification of $\M_0$ by nodes} \label{strat}

For every object $C \xrightarrow{\pi} T$ in $\M_0$, we are going to define the {\it relative singular locus of $C$}, which will be denoted as $C_{rs}$: roughly speaking it is the subfamily of $C$ whose geometric fibers have nodes. Its image to the base $T$ is a closed subscheme. Moreover it is possible to define on $T$ closed subschemes $\{ T^{\geq k} \to T\}_{k \geq 0}$ where for each integer $k \geq 0$ the fiber product $T^{\geq k} \times_T C$ is a family of rational curves with at least $k$ nodes.
In order to give a structure of closed subspaces to $C_{rs}$ and to the various subschemes $T^{\geq k}$, we will describe them through Fitting ideals.

\begin{defi} We construct the closed subspace $C_{rs}$ and the $T^{\geq k}$ locally in the Zariski topology, so we consider a family of rational nodal curve $C \xrightarrow{\pi} \spe A$ for some $\C-$algebra $A$.

Now we follow \cite{mumlect} Lecture 8.

As the map $\pi$ has relative dimension 1 the relative differential sheaf $\Omega_{C/T}$ has rank 1 as a sheaf over $C$. Further because $\pi$ is a map of finite presentation and $A$ is a Noetherian ring, we have that the sheaf $\Omega_{C/T}$ is coherent. For every point $p \in T$ we set
 \begin{eqnarray*}
e(p)=  \dim _{k(p)} (\Omega _{C/T,p} \otimes k(p))
\end{eqnarray*}
where $k(p)$ is the residual field of $p$.

Choose a basis $\{a_i\}_{i=1, \dots, e(p)}$ of $\Omega _{C/T,p} \otimes k(p)$, they extend to a generating system for $\Omega _{C/T,p}$. Furthermore we have an extension of this generating system to $\Omega_{C/T}$ restricted to an \'etale  neighborhood of $p$. So we have a map
$$
\Of_{C} ^{\oplus e(p)} \to \Omega_{C/T}
$$
which is surjective up to a restriction to a possibly smaller neighborhood of 
$p$. At last (after a possibly further restriction) we have the following exact sequence of sheaves
\begin{eqnarray*}
\Of_{C}^{\oplus f} \xrightarrow{X} \Of_{C} ^{\oplus e(p)} \to \Omega_{C/T} \to 0
\end{eqnarray*}
where $X$ is a suitable matrix $f \times e$ of local sections of $\Of_{C}$. Let us indicate with $F _i(\Omega_{C/T}) \subset \Of_C$ the ideal sheaf generated by rank $e(p) - i$ minors of $X$. These sheaves are known as  {\it Fitting sheaves}
and they don't depend on the choice of generators (see \cite{lan} XIX, Lemma 2.3).

We define the {\it relative singular locus} to be the closed algebraic subspace $C_{rs} \subseteq C$ associated with $F_1(\Omega_{C/T})$.
\end{defi}

\begin{oss}Let us fix a point $t \in T$, we have that $t$ belongs to $T^{\geq k}$ iff
$$
F_{k-1}(\pi_*(\Of_{C_{rs}})_t) = 0 \Longleftrightarrow \dim_{k(t)} (\pi_*(\Of_{C_{rs}})_t \otimes k(t)) \geq k
$$
so $T^{\geq k}$ is the subscheme whose geometric fibers have at least $k$ nodes.
\end{oss}

\begin{defi}Set
\begin{eqnarray} \label{stratk}
T^0 &:=& T - T^{\geq 1}\\
T^k &:=& T^{\geq k} - T^{\geq k + 1}
\end{eqnarray}
\end{defi}

\begin{oss}We have that for every $k \in \N$ the subscheme $T^k$ is locally closed in $T$, further, given a point $t \in T$ and set 
$$
k:= \dim_{k(t)} (\pi_*(\Of_{C_{rs}})_t \otimes k(t)) \geq 0
$$
we have that $t$ belongs to a unique $T^k$.

Consequently, given a curve $C \xrightarrow{\pi} T$ where $T$ is an affine Noetherian scheme, the family 
$$
\{ T^k | T^k \neq \emptyset \}
$$
defined above is  {\it a stratification for $T$}.
\end{oss}

Now let $\St$ be an Artin stack over $\cat$. A family of locally closed substacks
$\{ \St^{\alpha}\}_{\alpha \in \N}$ represents a stratification for $\St$ if, for all morphisms
$$
T \to \St
$$
where $T$ is a scheme, the family of locally closed subschemes
$$
\{ T^{\alpha} := \St^{\alpha} \times_{\St} T \neq \emptyset \}
$$
is such that every point $t \in T$ is in exactly one subscheme $T^{\alpha}$, that is to say that $\{ T^{\alpha}\}$ is a stratification for $T$ in the usual sense.

\begin{defi}

We define $\M_0 ^{\geq k}$ as the full subcategory of $\M_0$ whose objects $C \xrightarrow{\pi} T$ are such that $T^{\geq k} = T$.

We further define $\M_0 ^{k}$ to be the subcategory whose objects $C \xrightarrow{\pi} T$ are such that $T^{k} = T$.

\end{defi}

\begin{oss}
We have from definition that for every morphism $T \to \M_0$
$$
T^{\geq k}=\M^{\geq k}_0 \times_{\M_0} T
$$
consequently $\{ \M^{k}_0 \}_{k \in \N}$ is a stratification for $\M_0$.
\end{oss}

\begin{propos} \label{stratreg}
For each ${k \in \N}$ the morphism
$$
\M_0 ^{\geq k} \xrightarrow{} \M_0
$$
is a regular embedding of codimension $k$.
\end{propos}

\begin{dem}
Let us consider a smooth covering
$$
U \xrightarrow{f} \M_0
$$
and let $C \xrightarrow{\pi} U$ be the associated curve.
we have to prove that the morphism
$$
U^k= \M_0 ^k \times_{\M_0} U \to U
$$
is a regular embedding of codimension $k$ (we already know that $U^k$ is a closed embedding).

Fix a point $p \in U^k\subset U$ and let $A:=\widehat{\Of}_{U,p} ^{sh}$ be the completion of the strict henselisation of the local ring $\Of_{U,p}$. (see \cite{egaq} Definiton 18.8.7). Consider the following diagram
\begin{eqnarray*}
\xymatrix{
C_A \cart \ar[d]_{\pi} \ar[r] & C \ar[d]_{\pi}\\
\spe A \ar[r] & U
}
\end{eqnarray*}
As $p$ is a point of $U^k$ we have that $C_{A, rs}$ is the union of $k$ nodes $q_1, \dots, q_k$. For each $i=1, \dots, k$ we have
$$
\widehat{\Of}_{C_A, q_i} =A \llbracket x,y \rrbracket /(xy-f_i)
$$
with $f_1, \dots , f_k \in A$.
So we can write
$$
M:=\Of_{C_{A,rs}} = \prod_{i=1}^k (A/f_i).
$$
Let us consider $M$ as an $A-$module, we have the following exact sequence of $A-$module
$$
A^k \xrightarrow{D} A^k \xrightarrow{} M \to 0
$$
where $D$ is the diagonal matrix
\begin{equation}
\left(\begin{array}{ccc}
f_1 & \cdots & 0 \\
\vdots & \ddots & \vdots \\
0 & \cdots & f_k
\end{array}\right)
\end{equation}
So we have $F_{k-1}(M)=(f_1, \dots, f_k)$ and this means that
$$
B:=\widehat{\Of}^{sh}_{U^{\geq k},p} = A/(f_1, \dots, f_k).
$$
From deformation theory we have that $\{ f_1, \dots, f_k \}$ is a regular sequence, consequently the map $U^{\geq k} \to U$ is regular of codimension $k$ as claimed.
\end{dem}

\section{Combinatorical version of rational nodal curves.}
Now let us fix a useful notation. Given an algebraically closed field $\Omega$ and an isomorphism class of $C$ (still denoted with $C$) in $\M_0(\spk)$, we define the {\it dual graph of $C$}, denoted $\Gamma(C)$ or simply $\Gamma$, to be the graph which has as many vertices as the irreducible components of $C$ and two vertices are joined by an edge if and only if the two corresponding lines meet each other.

For example we have the following correspondence
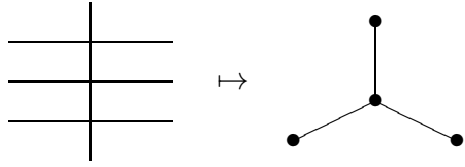
\begin{eqnarray*}
\vcenter{\trebg} & \; \mapsto \; & \vcenter{\treb}
\end{eqnarray*}
We can associate at least one curve with every tree by this map, but such a curve is not in general unique up to isomorphism; an example is the following tree
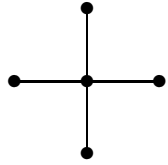
\begin{equation*}
\quattrob
\end{equation*}
In the following we think of a tree $\Gamma$ as a finite set of vertices with a connection law given by the set of pairs of vertices which corresponds to the edges.

\begin{defi}
Given a graph $\Gamma$ we call {\it multiplicity} of a vertex $P$ the number $e(P)$ of edges to which it belongs and we call $E(P)$ the set of  edges to which it belongs. Furthermore we call the {\it maximal multiplicity} of the graph $\Gamma$ the maximum of multiplicities of its vertices. We also call $\Delta_n$ the set of vertices with multiplicity $n$ and $\delta_n$ the cardinality of $\Delta_n$.
\end{defi}

\begin{oss} \label{isclass}
We can associate an unique isomorphism class of curves in $\M_0(\spk)$ to a given tree  $\Gamma$ if and only if the maximal multiplicity of $\Gamma$ is 3.
\end{oss}

Now let us fix a stratum $\M_0^k$, there are as many topological types of curves with $k$ nodes as trees with  $k+1$ vertices.

Purely topological arguments show the following

\begin{lemma} Given a curve $C \xrightarrow{\pi} T$ in $\M_0^{k}$ where $T$ is a connected scheme, the curves of the fibers are of the same topological type. 
\end{lemma}

So we can give the following:

\begin{defi}
For each tree $\Gamma$ with $k+1$ vertices we define $\M_0^{\Gamma}$ as the full subcategory of $\M_0 ^k$ whose objects are curves $C \xrightarrow{\pi} T$ such that on each connected component of $T$ we have curves of topological type $\Gamma$.
\end{defi}

For each $\Gamma$, $\M^{\Gamma}_0$ is an open (and closed) substack of $\M^{k}_0$ and we can consequently write
$$
\M_0 ^k = \coprod_{\Gamma} \M_0 ^{\Gamma}
$$
where $\Gamma$ varies among trees with $k+1$ vertices.

\section{Description of particular strata}\label{strata}

From now on we focus on graphs (and curves) with maximal multiplicity equal to 3. 

\begin{lemma}
Let $\Gamma$ be a graph with maximal multiplicity equal to 3 and $C$ an isomorphic class of curves of topological type $\Gamma$, then we have the equivalence
\begin{equation} \label{equiv}
\M_0 ^{\Gamma} \simeq \B \aut(C).
\end{equation}
\end{lemma}

\begin{dem}
From Remark (\ref{isclass}) we have that all curves of the same topological type $\Gamma$ are isomorphic.
\end{dem}

We have a canonical surjective morphism
$$
\text{Aut}(C) \xrightarrow{g} \text{Aut}(\Gamma)
$$
which sends each automorphism to the induced graph automorphism.

\begin{propos} \label{coord}
There is a (not canonical) section $s$ of $g$
\begin{equation*}
\xymatrix{
\text{Aut}(C) \ar@<-1ex>[r]_g &  \text{Aut}(\Gamma) \ar@<-1ex>[l]_s 
}
\end{equation*}
\end{propos}

\begin{dem}
Let us fix coordinates $[X,Y]$ on each component of $C$ such that
\begin{itemize}
\item on components with one node the point $[1,0]$ is the node,
\item on components with two nodes the points $[0,1]$ and $[1,0]$ are the nodes,
\item on components with three nodes the points $[0,1]$, $[1,1]$ and $[1,0]$ are the nodes.
\end{itemize}
Let us define on each component
\begin{eqnarray*}
0:=[0,1] \quad 1:=[1,1] \quad \infty:=[1,0]
\end{eqnarray*}

In order to describe the section $s$, let us notice that, given an element $h$ of $\aut(\Gamma)$, there exists a unique automorphism $\gamma$ of $C$ such that:
\begin{itemize}
\item $\gamma$ permutes components of $C$ by following the permutation of vertices given by $h$
\item on components with one node, $\gamma$ makes correspond the points $0,1$
\item on components with two nodes, $\gamma$ makes correspond the point $1$.
\end{itemize}
\end{dem}

\begin{propos}
We have
\begin{equation*}
\aut (C) \cong \aut(\Gamma) \ltimes \left( (\Gm)^{\Delta_2} \times E^{\Delta_1} \right).
\end{equation*}
where $E$ is the subgroup of $\Pro GL_2$ that fixes $\infty$ and $\Gm$ the multiplicative group of the base field.
\end{propos}

\begin{dem}
Let us consider the normal subgroup $g^{-1}(\id)$ of automorphisms of $C$ which do not permute components. It is the direct product of groups of automorphisms of each component that fixes nodes.

On components with one node the group of automorphisms is the subgroup $E$ of $\Pro GL_2$ that fixes $\infty$. $E$ can be described as the semidirect product of $\Gm$ and $\Ga$ (the additive group of the base field), moreover, having fixed coordinates, we can have an explicit split sequence
$$
\xymatrix{
0 \ar[r] & \Ga  \ar[r]^{\varphi} & E \ar@<-1ex>[r]_{\rho} & \Gm \ar@<-1ex>[l]_{\psi} \ar[r] & 1.
}
$$

On components of two nodes we have that the group of automorphisms is $\Gm$.

At last only identity fixes three points in $\Pro^1_\Omega$. So we can conclude that
\begin{equation*}
g^{-1}(\id)=E^{\Delta_1} \times \Gm^{\Delta_2}
\end{equation*}
and we have a (not canonical) injection
\begin{equation*}
E^{\Delta_1} \times \Gm^{\Delta_2} \to \text{Aut}(C)
\end{equation*}

Then $Aut(C)$ is the semi-direct product given by the exact sequence
\begin{equation} \label{seps}
\xymatrix{
1 \ar[r] & \Gm^{\Delta_2} \times E^{\Delta_1} \ar[r] & \text{Aut}(C) \ar@<-1ex>[r] &  Aut(\Gamma) \ar@<-1ex>[l] \ar[r] & 1
}
\end{equation}
and we write it as
\begin{equation*}
\text{Aut}(\Gamma) \ltimes (\Gm)^{\Delta_2} \times E^{\Delta_1}.
\end{equation*}
\end{dem}

So we can explicit
\begin{equation} \label{ps}
\M^{\Gamma}_0= \B \left( \text{Aut}(\Gamma) \ltimes (\Gm)^{\Delta_2} \times E^{\Delta_1} \right).
\end{equation}

\section{Dualizing and normal bundles} \label{normal}

In this section we briefly give the description of basic bundles over $\M_0$ and its strata $\M^k_0$. In particular we point out their restriction to $\M_0^{\leq 3}$.

{\bf The dualizing bundle.} On $\M_0$ we consider the universal curve $\cu \xrightarrow{\Pi} \M^{\Gamma}_0$.
defined in the following way:
\begin{defi}
Let $\cu$ be the fibered category on $\cat$ whose objects are families $C \xrightarrow{\pi} T$ of $\M_0$ equipped with a section $T \xrightarrow{s}C$ and whose arrows are arrows in $\M_0$ which commute with sections.
\end{defi}

Let $U \to \M_0$ be a smooth covering ($\M_0$ is an Artin stack) and let us consider the following cartesian diagram

\begin{equation*}
\xymatrix{
 \cu_U \cart \ar[d]^\Pi \ar[r] & \cu \ar[d]^\Pi\\
U \ar[r] & \M_0
}
\end{equation*}

\begin{defi}
We define the relative dualizing sheaf $\omega_0$ of $\cu \xrightarrow{\Pi} \M_0$ as the dualizing sheaf $\omega_{\cu_U/ U}$ on $\cu_U$.
\end{defi}

This is a good definition because for each curve $C \xrightarrow{\pi} T$ in $\M_0$ there is  the dualizing sheaf $\omega_{C/T}$ and its formation commutes with the base change.

For each curve $C \xrightarrow{\pi} T$ in $\M_0^{\leq k}$ we have the push forward $\pi_* \omega^{\vee}_{C/T}$. If we wish to have a well defined push forward of $\left ( \omega_{0} \right )^{\vee}$ we need to prove that for each curve in $\M_0^{\leq k}$ the sheaf $\pi_{*} \omega ^{\vee}_{C/T}$ is locally free of constant rank and it respects the base change. We can do it when $k$ is 3. In particular we have the following

\begin{propos} \label{dualsh}
Let $C \xrightarrow{\pi} T$ be a curve in $\M_0^{\leq 3}$. Then $\pi_{*} \omega_{C/T}^{\vee}$ is a locally free sheaf of rank 3 and its formation commutes with base change.
\end{propos}

\begin{dem}
First of all we prove that it is locally free and its formation commutes with base change.
It is enough to show that for every geometric point $t$ of $T$ 
$$
H^1(C_t, \omega^{\vee}_{C_t/t})=0.
$$
From Serre duality
$$
H^1(C_t, \omega^{\vee}_{C_t/t})=H^0(C_t, \omega^{\otimes 2}_{C_t/t})^{\vee}.
$$
When the fiber is isomorphic to $\Pro ^1$ we have 
$$
\omega^{\otimes 2}_{C_t/t} = (\Of(-2))^{\otimes 2} = \Of(-4)
$$
and we do not have global sections different from 0.

When $C_t$ is singular we have that the restriction of $\omega^{\otimes 2}_{C_t/t}$ to components with a node is $(\Of(-1))^{\otimes 2}=\Of(-2)$ and consequently the restriction of sections to these components must be 0. The restriction to components with two nodes is $\Of(0)$ (so restriction of sections must be constant) and restriction to components with three nodes is $\Of(2)$ (in this case the restriction of sections is a quadratic form on $\Pro^1$).

Given these conditions, global sections of $\omega^{\otimes 2}_{C_t/t}$ on curves with at most three nodes (as they have to agree on the nodes)
must be zero.

Similarly we verify that $h^0(C_t, \omega^{\vee}_{C_t/t}) = 3$ and conclude by noting that
$$
(\pi_*(\omega^{\vee}_{C/T}))_t = H^0 (C_t, \omega^{\vee}_{C_t/t}).
$$
\end{dem}

We notice that in $\M_0^{\leq 4}$ there are curves $C \xrightarrow{\pi} \spk$ for which there are global sections for $\omega^{\otimes 2}_{C/\spk}$. When we have a curve of topological type
$$
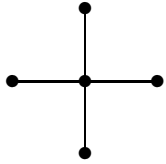
\quattrob
$$
the restriction of global sections on the central component are quartic forms on $\Pro^1_{\Omega}$ which vanishes on four points, thus we have
$$
H^0(C, \omega^{\otimes 2}_{C/\spk}) = \Omega.
$$ 

In conclusion we have a well defined bundle $\Pi_*\left ( \omega_{0} \right )^{\vee}$ on $\M^{\leq 3}_0$ and we refer to it as the dualizing bundle of $\M^{\leq 3}_0$.

\medskip

{\bf The normal bundles.} Given a tree $\Gamma$ we consider the (local) regular embedding (see Proposition \ref{stratreg})
$$
\M_0^\Gamma \xrightarrow{in} \M_0^{\leq \delta}
$$
where $\delta$ is the number of edges in $\Gamma$.
From \cite{kre} Section 5 we have that there exists a relative tangent bundle $\mathcal T_{in}$ on $\M_0^\Gamma$ which injects in $in^*(\mathcal T_{\M_0^{\leq \delta}})$.
\begin{defi}
Let us consider the following exact sequence of sheaves
$$
0 \to \mathcal T_{in} \to in^*(\mathcal T_{\M_0^{\leq \delta}}) \to in^*(\mathcal T_{\M_0^{\leq \delta}}) / \mathcal T_{in} \to 0.
$$
We define the normal bundle as the quotient sheaf on $\M_0^\Gamma$
$$
\Nf_{in}:=in^*(\mathcal T_{\M_0^{\leq \delta}}) / \mathcal T_{in}
$$
\end{defi}
When $\Gamma$ has maximal multiplicity at most 3, we can describe it as the quotient of the space of first order deformations by $\aut(C)$ in the following way.
Let us consider the irreducible components $ C_{\Gamma}$ of $C$. The space of first order deformations near a node $P$ where two curves $C_\alpha$ and $C_\beta$ meet is (see \cite{HM} p. 100)
$$
T_P(C_\alpha) \otimes T_P(C_\beta).
$$
Consequently we have the following:

\begin{lemma}
The space $N_\Gamma$ of first order deformations of $C$ is
$$
{\bigoplus_{P \in E(\Gamma)}} T_P(C_\alpha) \otimes T_P(C_\beta).
$$
\end{lemma}

\end{document}